\documentclass[12pt,reqno]{amsart}

\newtheorem{theorem}{Theorem}[section]
\newtheorem{lemma}[theorem]{Lemma}
\newtheorem{corollary}[theorem]{Corollary}
\theoremstyle{definition}   

\newtheorem{example}[theorem]{Example}

\theoremstyle{remark}

\numberwithin{equation}{section}

\usepackage{times}
\usepackage{enumerate}
\usepackage{graphicx,adjustbox}
\usepackage{hyperref}
\usepackage{amsmath,amssymb}

\usepackage{amscd}
\usepackage{graphicx}
\usepackage[all]{xy}

\title[Tangent Cones of Concatenated Numerical Semigroups]
{Tangent Cones of Concatenated Numerical Semigroups }
\author{
Ranjana Mehta
\and
Joydip Saha
\and
Indranath Sengupta
}
\date{}

\address{\small \rm  Department of Mathematics,
SRM University - AP, Amaravati, Andhra Pradesh 522502, India}
\email{ranjana.m@srmap.edu.in}

\address{\small \rm  Department of Mathematics, Barasat College, 1 Kalyani, Road, Barasat, West Bengal, PIN 700126 India.} 
\email{saha.joydip56@gmail.com}

\address{\small \rm  Discipline of Mathematics, IIT Gandhinagar, Palaj, Gandhinagar, 
Gujarat 382355, INDIA.}
\email{indranathsg@iitgn.ac.in}
\thanks{The third author is the corresponding author, 
supported by the grant CRG/2022/007047.}

\date{}



\allowdisplaybreaks

\begin{document}

\begin{abstract}
We study the tangent cone at the origin and the Hilbert series for a family of numerical semigroups generated by concatenation of arithmetic sequences. We prove that all the concatenation classes have Cohen-Macaulay tangent cones except the symmetric class, however, 
the symmetric class does satisfy Rossi's conjecture. 
\end{abstract}
\maketitle

\section{introduction}
Let $(R, \mathfrak{m}, k)$ be a one-dimensional Cohen-Macaulay local 
ring with associated graded ring 
$G = \oplus_{n\geq 0}\left(\frac{\mathfrak{m}^{n}}{\mathfrak{m}^{n+1}}\right)$. 
Let $H_{R}$ be the Hilbert function defined as 
$H_{R}(n) = H_{G}(n) = \dim_{k}\left(\frac{\mathfrak{m}^{n}}{\mathfrak{m}^{n+1}}\right)$. 
It is well-known 
that if $G$ is Cohen-Macaulay then $H_{R}$ is non-decreasing. Extensive study 
has been carried out by several authors on the Cohen-Macaulayness of $G$ and the 
properties of the Hilbert function. Rossi \cite{RossSurv}; Problem 4.9,  conjectured: 
\textit{If $R$ is a Gorenstein one-dimensional local ring, is it true that the 
Hilbert function of the ring is not decreasing?} Despite having affirmative 
answer for several classes (see \cite{ars1}, \cite{ars2}, \cite{ars3}, 
\cite{ars4}, \cite{JZ}, \cite{ot}), it was disproved in a recent work of 
Oneto, Strazzanti \& Tamone \cite{ost}. Most of the examples 
and counter examples belong to the class of numercial semigroup rings 
and often in embedding dimension $4$. 
\medskip

In \cite{mss3}, the authors introduced the notion of \textit{concatenation of two arithmetic sequences} to define a new family of numerical semigroups. This definition was largely inspired by the family of numerical semigroups in embedding dimension $4$, defined by Bresinsky in \cite{bre}. What turned out 
to be interesting is that the concatenation of arithmetic sequences, defined in \cite{mss3}, gives rise to numerical semigroups with different properties, for example, families with 
unique Ap\'{e}ry expansion, families which are symmetric and also families with unbounded 
Betti numbers. Patil-Tamone \cite{pt} proved Rossi's conjecture for certain 
classes of balanced numerical semigroup rings in embedding dimension $4$ and 
concatenation gives a generalization of the balanced numerical semigroups in arbitrary 
embedding dimension. 
In this paper, we have studied the Cohen-Macaulayness of the tangent cone at the 
origin and the Hilbert series of numerical semigroup rings for the concatenation 
class. The most interesting result proved is in section 3, which 
is a class of symmetric concatenated numerical semigroup 
(named as the symmetric class in \cite{mss3}), having a non Cohen-Macaulay tangent cone 
and a non-decreasing Hilbert function for all embedding dimensions $e\geq 4$; hence supports Rossi's conjecture. For other concatenation classes, we show that the tangent cone is Cohen-Macaulay and hence 
have a non-decreasing Hilbert function. The main novelty of this work lies in the technique. 
We have computed the complete Ap\'{e}ry table for the numerical semigroups, in all the three cases, to 
derive the tangent cone. Usually, it is not easy to give a complete description 
of the Ap\'{e}ry table. 

\section{Preliminaries and Notations}
A \textit{numerical semigroup} $\Gamma$ is a subset of the set of nonnegative integers $\mathbb{N}$, closed under addition, contains zero and generates $\mathbb{Z}$ as a group. We refer to \cite{rgs} for basic facts on numerical semigroups. 
Let $\Gamma$ be a numerical semigroup. It is true that 
(see \cite{rgs}) the set $\mathbb{N}\setminus \Gamma$ is 
finite and that the semigroup $\Gamma$ has a unique minimal system of generators 
$ a_{1} < \cdots < a_{e}$. The greatest integer not belonging to $\Gamma$ 
is called the \textit{Frobenius number} of $\Gamma$, denoted by $F(\Gamma)$. The integers 
$a_{1}$ and $e$ are known as the \textit{multiplicity} and the 
\textit{embedding dimension} of the semigroup $\Gamma$, usually 
denoted by $m(\Gamma)$ and $e(\Gamma)$ respectively. The 
\textit{Ap\'{e}ry set} of $\Gamma$ with respect to a non-zero $a\in \Gamma$ is 
defined to be the set $\rm{Ap}(\Gamma,a)=\{s\in \Gamma\mid s-a\notin \Gamma\}$.
Each element $x \in \Gamma$ can be written as $x = \sum_{i=1}^{e} a_{i}s_{i}$
for some non-negative integers $s_{i}$. The vector $\mathbf{s}=(s_{1},\ldots,s_{e})$ is called a factorization of $x$ and the
set of all factorizations of $x$ is denoted by $F(x)$, which is obviously a finite set. Let
$\mid \mathbf{s}\mid =\sum_{i=1}^{e} s_{i}$ denote the total order of $\mathbf{s}$. Then the maximum integer $n$ which is the total order of a vector in $F(x)$ is called the order of $x$ and is denoted by $\mathrm{ord}_{\Gamma}(x)$.
A vector $\mathbf{s} \in F(x)$ with $\mid\mathbf{s}\mid = \mathrm{ord}_{\Gamma}(x)$, is called a maximal factorization of $x$ and $x = \sum_{i=1}^{e} a_{i}s_{i}$ is called a maximal expression of $x$. For a vector a of non-negative
integers, we set $x(\mathbf{s}) = \sum_{i=1}^{e} a_{i}s_{i}$.

Given $0\neq x \in \Gamma$, the set of lengths of $x$ in $\Gamma$ is defined as
$$L(x) =\{\sum_{i=1}^{e}r_{i} | x = \sum_{i=1}^{e}r_{i}a_{i},r_i \geq 0\}$$

A subset $T \subset \Gamma$ is called homogeneous if either it is empty or $L(x)$
is singleton for all $0\neq x\in T$ . In other words, all expressions of elements in $T$ are
maximal. The numerical semigroup $\Gamma$ is called homogeneous, when the Ap\'{e}ry  set
$\mathrm{AP}(\Gamma, a_{1})$ is homogeneous.

Given integers $ a_{1} < \cdots < a_{e}$; the map 
$\nu : k[x_{1}, \ldots, x_{e}]\longrightarrow k[t]$ defined as 
$\nu(x_{i}) = t^{a_{i}}$, $1\leq i\leq e$, defines a parametrization 
for an affine monomial curve; the ideal $\ker(\nu)=\mathfrak{p}$ is called the 
defining ideal of the monomial curve defined by the parametrization 
$\nu(x_{i}) = t^{a_{i}}$, $1\leq i\leq e$. The defining ideal 
$\mathfrak{p}$ is a graded ideal with respect to the weighted gradation 
and therefore any two minimal generating sets of $\mathfrak{p}$ have the 
same cardinality.

Suppose $M=\Gamma\setminus\{0\}$ and for a positive integer $n$, we write 
$nM:=M+\cdots+M$ ($n$-copies). Let $r:=min\{r|(r+1)M=a_{1}+rM\}$, this $r$ is called reduction number. Let $\mathfrak{m}$ be the maximal ideal of the 
ring $k[[t^{a_{1}},\ldots t^{a_{e}}]]$. Then $(n+1)M=a+nM$ for all $n\geq r$ 
if and only if $r=r_{(t^{a_{1}})}(\mathfrak{m})$.

Let $\mathrm{Ap}(\Gamma,a_{1})=\{0,\omega_{1},\ldots,\omega_{a_{1}-1}\}$. Now for 
each $n\geq 1$, let us define $\mathrm{Ap}(nM)=\{\omega_{n,0},\ldots\omega_{n,a_{1}-1}\}$ 
inductively. We define $\omega_{1,0}=a_{1}$ and $\omega_{1,i}=\omega_{i}$, for $1\leq i\leq a_{1}-1$. 
Then $\mathrm{Ap}(M)=\{a_{1},\omega_{1},\ldots,\omega_{a_{1}-1}\}$. Now we define 
$\omega_{n+1,i}=\omega_{n,i}$, if $\omega_{n,i}\in (n+1)M $, and $\omega_{n+1,i}=\omega_{n,i}+a_{1}$, 
otherwise. We note that $\omega_{n+1,i}=\omega_{n,i}+a_{1} $ for all $0\leq i\leq a_{1}-1$ and $n\geq r_{(t^{a_{1}})}(\mathfrak{m})$. Then, the Ap\'{e}ry table $\mathrm{AT}(\Gamma,a_{1})$ of $\Gamma$ is a table 
of size $(r_{(t^{a_{1}})}(\mathfrak{m})+1)\times a_{1}$, whose $(0,t)$ entry is $\omega_{t}$, 
for $0\leq t\leq {a_{1}-1}$ (we take $\omega_{0}=0$), and the $(s,t)$ entry is $\omega_{st}$, 
for $1\leq s\leq r_{(t^{a_{1}})}(\mathfrak{m})$ and $ 0\leq t\leq {a_{1}-1}$.

We take some definitions from \cite{cz2}. Let $W =\{a_{0},\ldots,a_{n}\}$ be a set of integers. We call it a \textit{ladder} if $a_{0}\leq\ldots\leq a_{n}$. Given a ladder, we say that a subset $L=\{a_{i},\ldots,a_{i+k}\}$, with $k\geq 1$, is a \textit{landing} of length $k$ if $a_{i-1}<a_{i}=\cdots=a_{i+k}<a_{i+k+1}$ (where $a_{-1}= -\infty$ and $a_{n+1}=\infty$). In this case, $s(L)=i$ and $e(L)=i+k$. A landing $L$ is said to be a \textit{true landing} if $s(L)\geq 1$. Given two landings $L$ and $L^{'}$, we set $L<L^{'}$ if $s(L)<s(L^{'})$. Let $p(W)+1$ be the number of landings and assume that $L_{0}<\cdots<L_{p(W)}$ are the distinct landings. Then we define the following numbers:
$s_{j}(W)=s(L_{j})$, $e_{j}(W)=e(L_{j})$, for each $0\leq j\leq p(W)$;
$c_{j}(W)=s_{j}(W)-e_{j-1}(W)$, for each $0\leq j\leq p(W)$.

Suppose $\Gamma$ be a numerical semigroup minimally 
generated by $a_{1}<\cdots <a_{e}$ and $\mathfrak{m}_{\Gamma}$ be the maximal ideal of $k[[t^{a_{1}},\ldots t^{a_{e}}]]$. Let $r= r_{(t^{a_{1}})}(\mathfrak{m}_{\Gamma})$,  $M=\Gamma\setminus\{0\}$ and 
$\mathrm{Ap}(nM)=\{\omega_{n,0},\ldots\omega_{n,a_{1}-1}\}$ for $0\leq n \leq r$. For every $1\leq i\leq a_{1}-1$, consider the ladder of the values $W^{i}=\{\omega_{n,i}\}_{0\leq n\leq r}$ and define the following integers:
\begin{enumerate}[(i)]
\item $p_{i}=p(W^{i})$
\item $d_{i}=e_{p_{i}}(W^{i})$
\item $b_{j}^{i}=e_{j-1}(W^{i})$ and 
$c_{j}^{i}=c_{j}(W^{i})$, for $1\leq j\leq p_{i}$.
\end{enumerate}
\medskip

\begin{theorem}\textbf{(Cortadellas, Zarzuela.)}\label{tangentcone} With the above notations, 
$$G_{\mathfrak{m}_{\Gamma}}\cong F_{\Gamma}\oplus\displaystyle\bigoplus_{i=1}^{a_{1}-1}\left(F_{\Gamma}(-d_{i})\displaystyle \bigoplus_{j=1}^{p_{i}}\dfrac{F_{\Gamma}}{(({t^{a_{1}})^{*})^{c_{j}^{i}}}F_{\Gamma}}(-b_{j}^{i})\right),$$
where $G_{\mathfrak{m}_{\Gamma}}$ is the tangent cone of $\Gamma$ and $F_{\Gamma}=F((t^{a_{1}}))$ is the fiber cone.
\end{theorem}

\proof See Theorem 2.3 in \cite{cz2}.\qed
\medskip

\noindent\textbf{Concatenation.} Let $e\geq 4$. Let us consider the string of positive integers 
in arithmetic progression: 
$a<a+d<a+2d<\ldots<a+(n-1)d<b<b+d<\ldots<b+(m-1)d$, where $m,n\in \mathbb{N}$, $m+n=e$  
and $\gcd(a,d)=1$. Note that $a<a+d<a+2d<\ldots<a+(n-1)d$ and 
$b<b+d<\ldots<b+(m-1)d$ are both arithmetic sequences with the 
same common difference $d$. We further assume that this sequence minimally generates 
the numerical semigroup $\Gamma=\langle a, a+d, a+2d,\ldots, a+(n-1)d, b, b+d,\ldots, b+(m-1)d\rangle$. 
Then, $\Gamma$ is called the numerical semigroup generated 
by concatenation of two arithmetic sequences with the same common difference $d$. 

\section{Tangent cone of Symmetric Class of Numerical Semigroups}
Let us consider the symmetric class of numerical semigroup given in \cite{mss3}. Suppose $e \geq 4$ be an integer, $q$ a positive integer and $m = e+2q+1$. Let $d$ be a positive integer that satisfies $gcd(m, d)=1$. Let us define $\mathcal{S}=\{m,m+d,(q+1)m+(q+2)d,(q+1)m+(q+3)d,\ldots,(q+1)m+(q+e-1)d \}$ and $\Gamma_{(e,q,d)}(\mathcal{S})$ be the numerical semigroup generated by $\mathcal{S}$.

\begin{theorem}\label{aperyS}
 The Ap\'{e}ry set $\mathrm{Ap}(\Gamma_{(e,q,d)}(\mathcal{S}),m)$ for the numerical 
semigroup $\Gamma_{(e,q,d)}(\mathcal{S})$ with respect to the element $m$ is 
$\beta_{1}\cup\beta_{2}\cup\beta_{3}$, where
\begin{align*}
\beta_{1}& =\{k(m+d)\mid 0\leq k\leq q+1\},\\
\beta_{2}& =\{k(m+d)+(q+1)m+(q+e-1)d \mid 0\leq k\leq q+1\},\\
\beta_{3}& =\{(q+1)m+(q+i)d\mid 2\leq i\leq e-2\}.
\end{align*} 
\end{theorem}

\proof See Theorem $4.2$ \cite{mss3}. \qed

\begin{theorem}\label{uniqueS}
The following statements are true for the 
numerical semigroup $\Gamma_{(e,q,d)}(\mathcal{S})$. 
\begin{enumerate}
\item[(i)] Every element in the set $\mathrm{Ap}(\Gamma_{(e,q,d)}(\mathcal{S}),m)\setminus\{(q+1)(m+d)+(q+1)m+(q+e-1)d\}$ has unique expansion, hence the given set is homogeneous.
\smallskip

\item[(ii)] $\mathrm{ord}((q+1)(m+d)+(q+1)m+(q+e-1)d))=q+2$.
\smallskip

\item[(iii)]$\mathrm{ord}(k(m+d)+jm)=k+j$, \, for \, all \, $0\leq k\leq q+1$, \, $j\geq 1$.
\smallskip

\item[(iv)] For \, $2\leq i\leq e-2$,
$$\mathrm{ord}( (q+1)m+(q+i)d+jm) = 
\begin{cases} j+1, & \quad\mathrm{if}\quad 0\leq j<i-1;\\
q+j+1, & \quad\mathrm{if}\quad j\geq i-1.
\end{cases}$$

\item[(v)] For \, $0\leq k\leq q+1$, 
$$\mathrm{ord}( k(m+d)+(q+1)m+(q+e-1)d+jm) =                    
\begin{cases}
k+1+j, & \quad\mathrm{if}\quad j<e-2;\\
k+q+j+1, & \quad\mathrm{if}\quad j\geq e-2.
\end{cases}$$

\end{enumerate}
\end{theorem}
\proof  We write $n_{1}=m, n_{2}=m+d$ and $n_{i}=(q+1)m+(q+i-1)d$, $ 3\leq i\leq e$.
\smallskip

\begin{enumerate}
\item[(i)] We have proved 
that each element in $\mathrm{Ap}(\Gamma_{(e,q,d)}(\mathcal{S}),m)$ except the maximal element has 
unique expansion in Theorem 4.5 \cite{mss3}, i.e., the set 
$$\mathrm{Ap}(\Gamma_{(e,q,d)}(\mathcal{S}),m)\setminus\{(q+1)(m+d)+(q+1)m+(q+e-1)d\}$$ is homogeneous.
\medskip
 
\item[(ii)] The maximal element of $\mathrm{Ap}(\Gamma_{(e,q,d)}(\mathcal{S}),m)$ is $(q+1) n_{2}+n_{e}$. Suppose $(q+1) n_{2}+n_{e}=\displaystyle\sum_{i=2}^{e}k_{i}n_{i}$, then we need to show that $\displaystyle\sum_{i=2}^{e}k_{i}\leq (q+2)$. We consider the following cases:
\medskip

\noindent\textbf{Case (a).}\, If $k_{e}=0$, we have
\begin{align*}
&(q+1)(m+d)+(q+1)m+(q+e-1)d\\
&=2(q+1)m+(2q+e)d\\
&=(k_{2}+\sum_{l=3}^{e-1}k_{l}(q+1))m+(k_{2}+\sum_{l=3}^{e-1}k_{l}(q+l-1))d
\end{align*}
Suppose $k_{2}=0$.
If $\displaystyle\sum_{i=3}^{e-1}k_{i}< 2$ then $\displaystyle\sum_{i=2}^{e}k_{i}< q+2$ and  we are done. If $\displaystyle\sum_{i=2}^{e-1}k_{i}\geq 2$ then we have,
\begin{eqnarray*}
 (\sum_{l=3}^{e-1}k_{l}-2)(q+1)m =((2q+e)-(\sum_{l=3}^{e-1}k_{l}(q+l-1))d.
\end{eqnarray*}
 
\noindent The R.H.S is at most $e-q-6$. We have 
$gcd(m,d)=1$, therefore $m|(e-q-6)$. However, 
$m=e+2q+1$, which is contradictory to $m|e-q-6$.
\medskip

If $k_{2}\neq 0$, then it follows that 
$0< k_{2}\leq q$ 
from the equation $((q+1)-k_{2}) n_{2}+n_{e}=\displaystyle\sum_{i=3}^{e-1}k_{i}n_{i}$. Therefore, if 
$\displaystyle\sum_{i=2}^{e-1}k_{i}\leq 2$, then 
$\displaystyle\sum_{i=3}^{e}k_{i}\leq (q+2)$.
\medskip

\noindent\textbf{Case (b).}\, If $k_{e}>0$, then $(q+1) n_{2}+n_{e}=\displaystyle\sum_{i=2}^{e-1}k_{i}n_{i}+k_{e}n_{e}$ implies
$$
(q+1)n_{2}-(k_{e}-1)n_{e}=(2-k_{e})qd+d-(k_{e}-1)(e-1)d
=\displaystyle\sum_{i=2}^{e-1}k_{i}n_{i}.$$

\noindent If $k_{e}\geq 2$, then 
$(2-k_{e})qd+d-(k_{e}-1)(e-1)d<0$, whereas,  
$\displaystyle\sum_{i=2}^{e-1}k_{i}n_{i}>0$, a contradiction. Therefore, $k_{e}=1$ and
$(q+1) n_{2}=\displaystyle\sum_{i=2}^{e-1}k_{i}n_{i}$. 
Since $n_{2}\leq n_{i}$, for $2\leq i\leq e$, we have $\displaystyle\sum_{i=2}^{e-1}k_{i}n_{2}\leq \displaystyle\sum_{i=2}^{e-1}k_{i}n_{i}=(q+1) n_{2}$. Hence 
$\displaystyle\sum_{i=2}^{e-1}k_{i}\leq (q+1)$. 
Since $k_{e}=1$, therefore, $\displaystyle\sum_{i=2}^{e}k_{i}\leq (q+2)$. Hence order of $(q+1) n_{2}+n_{e}$ is 
$(q+2)$.
\medskip

\item[(iii)] Next we wish to show that for $0\leq k\leq q+1$ and for all $j\geq 0$, $\mathrm{ord}(k(m+d)+jm)=k+j$. Suppose 
\begin{equation}\label{eq41}
k(m+d)+jm=\sum_{i=1}^{e} a_{i}n_{i}
\end{equation}
If $a_{1}\geq j$, then 
\begin{equation}
k(m+d)=(a_{1}-j)m+\sum_{i=2}^{e} a_{i}n_{i}
\end{equation}
Since $k(m+d)$ has unique expression, we have 
$a_{1}=j$, $a_{2}=k$ and $a_{i}=0$, for $3\leq i\leq e$ 
and we are done. 
\medskip
 
\noindent If $1\leq a_{1}<j$, then 
\begin{equation}
k(m+d)+(j-a_{1})m=\sum_{i=2}^{e} a_{i}n_{i},
\end{equation}
and we are done by induction.
\medskip

\noindent If $a_{1}=0$, then rearranging equation \ref{eq41}, we get
\begin{equation}
kd=\left(a_{2}+(q+1)(\sum_{i=3}^{e} a_{i})-k-j\right)m+ \left(a_{2}+(\sum_{i=3}^{e} a_{i}(q+i-1))\right)d
\end{equation}

We have $\gcd(m,d)=1$, therefore $a_{2}+(q+1)(\sum_{i=3}^{e} a_{i})-k-j=\ell d$. If $\ell<0$, then 
$\sum_{i=2}^{e} a_{i}<a_{2}+(q+1)(\sum_{i=3}^{e} a_{i})<k+j$ and we are done. If $\ell\geq 0$, then  
\begin{align*}
\sum_{i=2}^{e} a_{i}n_{i} &\geq (a_{2}+(q+1)(\sum_{i=3}^{e} a_{i}))(m+d)\\
&=(k+j+\ell d)(m+d) \, > \, k(m+d)+jm
\end{align*}
and we are done.
\medskip

\item[(iv)]  Suppose $n_{i+1}+jn_{1}=(q+1)m+(q+i)d+jm=\sum_{t=1}^{e}a_{t}n_{t}$. If $a_{1}\geq j$, then $n_{i+1}=(a_{1}-1)n_{1}+\sum_{t=2}^{e}a_{t}n_{t}$. The integer 
$n_{i+1}$ belongs to  the minimal generating set of $\Gamma_{(e,q,d)}(\mathcal{S})$, therefore 
$a_{i+1}=1, a_{1}=1,a_{t}=0$ for $t\neq 1,i+1$, and we are done.
\smallskip

If $1\leq a_{1}<j$, then $n_{i+1}+(j-a_{1})m=\sum_{t=2}^{e}a_{t}n_{t}$, and by induction we have $a_{2}+\cdots+a_{e}\leq 1+j-a_{1}$, i.e, $\sum_{t=1}^{e}a_{t}\leq j+1 $. Therefore, 
we assume that $a_{1}=0$. 
\smallskip

Next, if $a_{i+1}\geq 1$, then 
$$jm=\displaystyle\sum_{t=2,t\neq i+1}^{e}a_{t}n_{t}+(a_{i+1}-1)n_{i+1} \, > \, (\sum_{t=2}^{e}a_{t}-1)m,$$
 which shows that $j+1>\sum_{t=2}^{e}a_{t}$ and we are done. Therefore we further assume that $a_{i+1}=0$.
 \smallskip
 
 Next, if $a_{k}\geq 1 $ for some $i+2\leq k\leq e$, 
 and $j+1<\displaystyle\sum_{t=2,t\neq i+1}^{e}a_{t} $, then
 \begin{align*}
 jm &<\left(\displaystyle\sum_{t=2,t\notin\{i+1,k\}}^{e}a_{t}\right)+(a_{k}-1)m\\
 &<\left(\displaystyle\sum_{t=2,t\notin\{i+1,k\}}^{e}a_{t}n_{t}\right)+(a_{k}-1)n_{k},
 \end{align*}
which gives a contradiction. Therefore, 
$j+1\geq \sum_{t=2}^{e}a_{t}$ and we are done.
 \smallskip
 
 Hence, we may assume that $ a_{k}=0$, 
 for $i+2\leq k\leq e$, and 
 \begin{equation}\label{eq45}
 n_{i+1}+jn_{1}=\sum_{t=2}^{i}a_{t}n_{t}.
 \end{equation}
Rearranging terms of the equation \ref{eq45}, we get 
$$\left((q+i)-a_{2}-\sum_{t=3}^{i}(q+t-1)a_{t}\right)d=\left(a_{2}+(\sum_{t=3}^{i}a_{t})-1)(q+1)-j\right)m.$$
We have $\gcd(m,d)=1$, therefore, 
\begin{eqnarray}\label{eq46}
(q+i)-a_{2}-\sum_{t=3}^{i}(q+t-1)a_{t}=lm\\
a_{2}+(\sum_{t=3}^{i}a_{t})-1)(q+1)-j=ld
\end{eqnarray}
If $l\geq 1$, then $(q+i)<m$ \, for \, $2\leq i\leq e-2$, 
as $m=e+2q+1$. This contradicts equation \ref{eq46}, 
therefore $l\leq 0$. 
\smallskip
 
\noindent\textbf{Subcase 1.} If $a_{3}=\cdots=a_{i}=0$, then $a_{2}-(q+1)-j\leq 0$ i.e, $a_{2}\leq q+1+j$. Again $(q+i)\leq a_{2}$, implies that $q+i\leq q+1+j$ i.e, $j\geq i-1$.
\smallskip

\noindent\textbf{Subcase 2.} If $\sum_{t=3}^{i}a_{t}\geq 1$,  then $\sum_{t=2}^{i}a_{t}-1\leq a_{2}+(\sum_{t=3}^{i}a_{t}-1)(q+1)=j $, and we are done. Note that if $j\geq i-1$, then $n_{i+1}+jn_{1}=(q+i)n_{2}+(j-i-1)n_{1}$, so for a maximal presentation of $n_{i+1}+jn_{1} $ we have $\sum_{t=2}^{i}a_{t}\geq q+j+1>j+1$. Hence, for this subcase, we have $j<i-1$. 
\medskip

\item[(v)] Proof is similar as in (iv). \qed
\end{enumerate}

\begin{theorem}\label{aperytable(s)}
The Ap\'{e}ry table $\mathrm{AT}(\Gamma_{(e,q,d)}(\mathcal{S}), m)_{m\times m}$ of  $\Gamma_{(e,q,d)}(\mathcal{S})$ is of the form 
$\begin{pmatrix}
\mathcal{S}_{1}& \mathcal{S}_{2}& \mathcal{S}_{3}
\end{pmatrix}$, where $\mathcal{S}_{1}$, $\mathcal{S}_{2}$, $\mathcal{S}_{3}$ 
are described below:

$\bullet$ $\mathcal{S}_{1}=(s_{jk})_{m\times (q+2)}$, with
\begin{align*}
s_{1k}&=(k-1)(m+d),\, 1\leq k\leq (q+2),\\
s_{j1}&=(j-1)m_{0},\, 2\leq j\leq m,\\
s_{jk}&=(k-1)(m+d)+\delta_{jk}m,\, 2\leq k\leq (q+2),\, 2\leq j\leq m,
\end{align*}
such that \, 
$\delta_{jk}= \begin{cases}
0, & \quad \rm{if} \quad j<k;\\
j-k, & \quad \rm{if} \quad j\geq k, 2\leq k\leq (q+2), 2\leq j\leq m.
\end{cases}$
\smallskip

$\bullet$ $\mathcal{S}_{2}=(s_{jk})_{m\times (q+2)}$, with\\ 
$$s_{1k}  = (k-1)(m+d)+(q+1)m+(q+e-1)d, \quad \rm{if} \quad 1\leq k\leq (q+2);$$
$s_{jk} = 
\begin{cases}
s_{1k}, & \rm{if}\quad 1\leq k\leq (q+2),\,1\leq j\leq k+1;\\
s_{1k}+(j-k-1)m, & \rm{if} \quad 1\leq k\leq (q+2), k+2\leq j\leq e+k-2;\\
s_{1k}+(e-2)m, & \rm{if} \quad 1\leq k\leq (q+2), e+k-1\leq j\leq e+q+k-1;\\
s_{1k}+(j-q-k-1)m, & \rm{if} \quad 1\leq k\leq (q+2), e+q+k\leq j\leq m.
\end{cases}$
\smallskip

$\bullet$ $\mathcal{S}_{3}=(s_{jk})_{m\times (e-3)}$, with
$$s_{1k}=(q+1)m+(q+k+1)d,\, 1\leq k\leq e-3;$$
$s_{jk} = 
\begin{cases}s_{1k}+(j-2)m, & \quad \rm{if} \quad 1\leq k\leq (e-3),\, 2\leq j\leq k+1;\\
s_{1k}+km, & \quad \rm{if} \quad 1\leq k\leq (e-3), k+2\leq j\leq q+k+2;\\
s_{1k}+(j-q-2)m, & \quad \rm{if} \quad 1\leq k\leq (e-3), q+k+3\leq j\leq m.\\
\end{cases}$
\end{theorem}

\proof Follows from Theorem \ref{uniqueS}.\qed

\begin{corollary}\label{exptangent} Let $I=\langle t^{m}\rangle$. The tangent cone $G_{\mathfrak{m}}(\Gamma_{(e,q,d)}(\mathcal{S}))$ of $\Gamma_{(e,q,d)}(\mathcal{S})$ is a free $F(I)$-module. Moreover 
\begin{align*}
G_{\mathfrak{m}}(\Gamma_{(e,q,d)}(\mathcal{S}))=& \displaystyle\bigoplus F(I)\bigoplus_{i=1}^{q+1}F(I)(-i)\\
&\bigoplus_{i=1}^{q+2}\left(F(I)(-(e+q+i-2))\oplus\dfrac{F(I)}{(t^{m})^{e-2} F(I)}(-i)\right)\\
&\bigoplus_{i=1}^{e-3}\left(F(I)(-(q+i+1))\oplus\dfrac{F(I)}{(t^{m})^{i} F(I)}(-1)\right)
\end{align*} 
\end{corollary} 

\proof Proof follows from the theorems \ref{tangentcone}, \ref{unique4}, \ref{aperytable4} and Lemma \ref{order4}.\qed

\begin{corollary}Tangent cone $G_{\mathfrak{m}}(\Gamma_{(e,q,d)}(\mathcal{S}))$ is not Cohen-Macaulay.
\end{corollary}
\proof It follows from the fact that  $G_{\mathfrak{m}}(\Gamma_{(e,q,d)}(\mathcal{S}))$ is not a free $F(I)$-module (see section 4 in \cite{cz2}). \qed

\begin{theorem}\label{rossi} The Hilbert function $H_{\Gamma_{(e,q,d)}(\mathcal{S})}(n)$ of $\Gamma_{(e,q,d)}(\mathcal{S})$ is non decreasing in each embedding dimension $e\geq 4$, hence supports Rossi's conjecture. 
    \end{theorem}
    
\proof As given in section $2$\cite{bjz}, we have $H_{\Gamma_{(e,q,d)}(\mathcal{S})}(n)=\# (nM\setminus (n+1)M)$, the Hilbert function of $k[[{\Gamma_{(e,q,d)}(\mathcal{S})}]]$, where $M={\Gamma_{(e,q,d)}(\mathcal{S})}\setminus \{0\}$. From the theorem \ref{aperytable(s)}, we get, $\#(\mathrm{Ap}(\Gamma_{(e,q,d)}(\mathcal{S}),m)\setminus M)=1$ and
\begin{eqnarray*}
\# (nM\setminus (n+1)M) &=e+(n-1)\,\,\, \mathrm{for}\,\, 1\leq n\leq q+1 \\
&=e+q+1 \,\,\, \mathrm{for}\,\, q+2\leq n\leq e+q \\
&=n+1\,\,\, \mathrm{for}\,\, e+q+1\leq n\leq m-2 \\
\end{eqnarray*}
Hence, the Hilbert function of $k[[{\Gamma_{(e,q,d)}(\mathcal{S})}]]$ takes the values  
$$\{1,e,e+1,\ldots,e+q,e+q+1,e+q+1,\ldots(e-1)\,\rm{times},e+q+2,\ldots m-1,\rightarrow\}\qed$$ 
\begin{corollary}
Hilbert series of $\Gamma(e,q,d)(S)$ is
\begin{eqnarray*}
 H_{\Gamma(e,q,d)}(S)t & = & \sum_{i=1}^{q+1}(e+i-1)t^{i}+ \sum_{i = q+2}^{e+q}(e+q+1)t^{i}+\sum_{i = e+q+1}^{m-2}(n+1)t^{i}+\sum_{i = m-1}^{\infty}mt^{i} \\
 & = & \sum_{i=1}^{q+1}(e+i-1)t^{i}+ \sum_{i = q+2}^{e+q}(e+q+1)t^{i}+\sum_{i = e+q+1}^{m-2}(n+1)t^{i}+ \frac{mt^{m-1}}{(1-t)}\\
 & = & \frac{f(t)}{(1-t)}
 \end{eqnarray*}
 Where,
\begin{align*}
f(t)&=\displaystyle\sum_{i=1}^{q+1}(e+i-1)t^{i}(1-t)+ \displaystyle\sum_{i = q+2}^{e+q}(e+q+1)t^{i}(1-t)\\
&+\displaystyle\sum_{i=e+q+1}^{m-2}(n+1)t^{i}(1-t)+mt^{m-1}
\end{align*}
\end{corollary}
\proof Follows from the theorem \ref{rossi}.\qed
\begin{example}
Let $ S=\langle 9,11,35,37\rangle$,Here $e=4,q=2, m=9$ then 
$$\mathrm{AT}((S,9))_{9\times 9}=
\begin{tabular}{|c|c|c|c|c|c|c|c|c|c|} % centered columns (4 columns)
\hline %inserts double horizontal lines
$\mathrm{AP}(S)$ & $0$ & $11$ & $22$ & $33$ & $37$ & $48$ & $59$ & $70$ & $35$ \\ [0.5ex] 
\hline %inserts double horizontal lines
$M$ & $9$ & $11$ & $22$ & $33$ & $37$ & $48$ & $59$ & $70$ & $35$ \\ [0.5ex] 
\hline % inserts single horizontal line
$2M$ & $18$ & $20$  & $22$ & $33$ & $46$ & $48$ &$59$ & $70$ & $44$ \\ [0.5ex] 
\hline
$3M$ & $27$ & $29$ & $31$ & $33$ & $55$ & $57$ & $59$ & $70$ & $44$ \\ [0.5ex] 
\hline
$4M$ & $36$ & $38$ & $40$ & $42$ & $55$ & $66$ & $68$ & $70$ & $44$ \\ [0.5ex]
\hline
$5M$ & $45$ & $47$ & $49$ & $51$ & $55$ & $66$ & $77$ & $79$ & $53$ \\ [0.5ex]
\hline
$6M$ & $54$ & $56$ & $58$ & $60$ & $64$ & $66$ & $77$ & $88$ & $62$ \\ [0.5ex]
\hline
$7M$ & $63$ & $65$ & $67$ & $69$ & $73$ & $75$ & $77$ & $88$ & $71$ \\ [0.5ex]
\hline
$8M$ & $72$ & $74$ & $76$ & $78$ & $82$ & $84$ & $86$ & $88$ & $80$ \\ [0.5ex]
\hline
%$\Gamma_{4}(M)(t_{k})$  & $3$ & $1$ \\\hline % inserting body of the table
\end{tabular}$$ 
% is used to refer this table in the text

\noindent Therefore, 
\begin{itemize}
\item $M\setminus2M=\{9,11,37,35\}$,  $\#(M\setminus2M)=4$
\item $2M\setminus3M=\{18,20,22,46,48\}$,   $\#( 2M\setminus3M)=5$
\item $3M\setminus4M=\{27, 29, 31, 33 , 57 ,59\}$, $\#( 3M\setminus 4M)=6$
\item $4M\setminus 5M=\{36, 38,  40, 42, 68, 70, 44\}$, $\#( 4M\setminus 5M)=7$
\item $5M\setminus 6M=\{45, 47, 49, 51, 55, 79, 53\}$,   $\#( 2M\setminus3M)=7$
\item $6M\setminus 7M=\{54, 56, 58, 60, 64, 66, 62\}$, $\#( 6M\setminus 7M)=7$ 
\item $7M\setminus 8M=\{63, 65, 67, 69, 73, 75, 77, 71\}$, $\#( 7M\setminus 8M)=8$
\end{itemize}
\smallskip

\noindent The Hilbert Function takes the values $ \{ 1,4,5,6,7,7,7,8 \rightarrow\}$, in a 
non-decreasing pattern.\\

\noindent The Hilbert Series of the numerical semigroup $S$ is as follows
\begin{align*}
H_{\Gamma(4,2,2)}(S)t&=1+4t+5t^2+6t^3+7t^4+7t^5+7t^6+\displaystyle\sum_{i=7}^{\infty} 8t^{i}\\
&=1+4t+5t^2+6t^3+7t^4+7t^5+7t^6+\displaystyle\frac{8t^7}{1-t}\\
&=\displaystyle\frac{1+3t+t^{2}+t^{3}+t^{4}+t^{7}}{1-t}\\
\end{align*}
\end{example}

\section{Tangent cone of Almost Maximal Class of Numerical Semigroups}
In this section we consider almost maximal class of numerical semigroup that is given in \cite{mss3}. Let $e\geq 4$ be an integer; $a=e+1$, $b > a+(e-3)d$, $\gcd(a,d)=1$ 
and $d\nmid(b-a)$. Let $M=\{a, a+d, a+2d,\ldots,a+(e-3)d, b,b+d\}$ and 
we assume that the set forms a minimal generating set for the numerical 
semigroup $\Gamma_{e}(M)$, generated by the set $M$. We once again 
recall that definition of concatenation includes the minimality of 
the sequence of integers generating the numerical semigroup $\Gamma_{e}(M)$. 
For example, for the choice of $e=d=4, a=5, b=10$, we get the 
concatenated sequence $5, 9, 10, 14$, which is not a good example due to 
its non-minimality. Let us write $d\equiv i(\mathrm{mod}\, a)$. 

\begin{theorem}\label{apery} 
\begin{enumerate}
\item If $e=4,$ then 
\begin{enumerate}[(i)]
\item $\mathrm{Ap}(\Gamma_{4}(M),5)= \{0, 5+d, b, b+d, 2b \} $, \, if \, $b\equiv 2i (\mathrm{mod}\,5)$;

\item $\mathrm{Ap}(\Gamma_{4}(M),5)= \{0, 5+d, b, b+d,2(5+d)\}$, \, if \, $b\equiv 3i (\mathrm{mod}\, 5)$.
\end{enumerate}
\medskip

\item If $e\geq 5,$ then 
\begin{enumerate}[(i)]
\item $\mathrm{Ap}(\Gamma_{e}(M),a)=\{0, a+d, \ldots, a+(a-4)d, b, b+d, b+a+2d \} $, \, if \,  
$b\equiv(a-3)i (\mathrm{mod}\, a)$;

\item $\mathrm{Ap}(\Gamma_{e}(M),a)=\{0, a+d, \ldots, a+(a-4)d, b, b+d,2a+(a-3)d\}$, \, if \,  
$b\equiv(a-2)i (\mathrm{mod}\,a)$.
\end{enumerate}
\end{enumerate}
\end{theorem}

\proof See Theorem 5.3\cite{mss3}. \qed

\begin{theorem}\label{orderam}
\begin{enumerate}
\item If $e=4$, then for $j\geq 0$, we have
\begin{enumerate}[(i)]
\item $\mathrm{ord}((5+d)+5j)=j+1$;
\item $\mathrm{ord}(b+kd+5j)=j+1,\, 0\leq k\leq 1$;
\item $\mathrm{ord}((2b)+5j)=j+2,\, if \, b\equiv 2i(\mbox{mod}\,5)$;
\item $\mathrm{ord}(2(5+d)+5j)=j+2,\, if \, b\equiv 3i(\mathrm{mod}\, 5)$.
\end{enumerate}
\medskip

\item If $e\geq 5$, then for $j\geq 0$, we have
\begin{enumerate}[(i)]
\item $\mathrm{ord}((a+kd)+ja)=j+1,\, 0\leq k\leq (a-4)$;
\item $\mathrm{ord}((b+kd)+ja)=j+1,\, 0\leq k\leq 1$;
\item $\mathrm{ord}((b+a+2d)+ja)=j+2,\, if \, b\equiv(a-3)i(\mbox{mod}\, a)$;
\item $\mathrm{ord}(2a+(a-3)d+ja)=j+2,\, if \, b\equiv(a-2)i(\mathrm{mod}\,a)$. 
\end{enumerate}
\end{enumerate}
\end{theorem}

\proof From the proof of Theorem 5.3 in \cite{mss3}, 
we can prove the following results which give the proofs 
of the statements stated above:
\begin{enumerate}
\item Let $b\equiv 2i (\mathrm{mod}\,5)$ and 
\begin{equation}\label{eq5.1}
 2b=c_{1}(5+d)+c_{2}b+c_{3}(b+d).
 \end{equation}
From this equation we get $0\leq c_{2}\leq 1$, 
$0\leq c_{3}\leq 1$. We can not have 
$c_{2}=c_{3}=1$ because $c_{1}(5+d)+d>0$. 
Hence, either $c_{2}=0$ or $c_{3}=0$.
\smallskip
 
\noindent Taking modulo $5$ on each side of the equation 
$\ref{eq5.1}$, we get $4\equiv (c_{2}+3c_{3})(\mathrm{mod}\, 5)$. Therefore, $c_{2}+3c_{3}=4+5t$, where $t\geq 0$, which 
is not possible as $0\leq c_{2}\leq 1$ and $0\leq c_{3}\leq 1$ and one of them must be $0$. Hence, the element $2b\in \mathrm{Ap}(\Gamma_{4}(M),5)\setminus M$, is uniquely expressed and the order of $2b$ is $2$. 
\smallskip
 
Similarly we can show that if $b\equiv 3i(\mathrm{mod}\, 5)$, the element $2(5+d)\in \mathrm{Ap}(\Gamma_{4}(M),5)\setminus M$, is uniquely expressed and the order of $5+2d$ is $2$.
\medskip

\item Let $b\equiv(a-3)i(\mathrm{mod}\, a)$ and 
\begin{equation}\label{eq5.2}
a+b+2d=\displaystyle\sum_{k=1}^{a-4}c_{k}(a+kd)+c_{a-3}b+c_{a-2}(b+d).
\end{equation}

If $c_{a-2}\geq 2$, then $\displaystyle\sum_{k=1}^{a-4}c_{k}(a+kd)+c_{a-3}b+2(b+d)> a+b+2d$, which is a contradiction. 
Hence $0\leq c_{a-2}\leq 1$.

If $c_{a-2}=1$, then $a+d=\displaystyle\sum_{k=1}^{a-4}c_{k}(a+kd)+c_{a-3}b$, by equation \ref{eq5.2}. Hence, 
$c_{1}=1$ and $a+b+2d=(a+d)+(b+d)$.

If $c_{a-2}=0$, then $a+b+2d=\displaystyle\sum_{k=1}^{a-4}c_{k}(a+kd)+c_{a-3}b$, by equation \ref{eq5.2}. We have  
$b>a+2d$, hence $0\leq c_{a-3}\leq 1$.

\begin{enumerate}[(i)]  
\item If $c_{a-3}=1$, then $a+2d=\displaystyle\sum_{k=1}^{a-4}c_{k}(a+kd)$. It is possible only when $c_{2}=1$ and $a+b+2d=(a+2d)+b$.

\item If $c_{a-3}=0$, then $a+b+2d=\displaystyle\sum_{k=1}^{a-4}c_{k}(a+kd)$. If $c_{1}=1$ then $b+d=\displaystyle\sum_{k=2}^{a-4}c_{k}(a+kd)$, which is a contradiction to the minimality of the generating set $M$.

Hence, $c_{1}=c_{a-3}=c_{a-2}=0$. Substituting in the equation \ref{eq5.2} we get, 
 \begin{equation}\label{eq5.3}
a+b+2d=\displaystyle\sum_{k=2}^{a-4}c_{k}(a+kd).
\end{equation}
Taking modulo $a$ both the sides of the equation \ref{eq5.3}, we get 
\begin{equation}\label{eq5.4}
(a-1)\equiv \displaystyle\sum_{k=2}^{a-4}c_{k}k(\mathrm{mod}\, a)
\end{equation}

Hence $\displaystyle\sum_{k=2}^{a-4}c_{k}k=ta+(a-1)$, where $t\geq 0$. Substituting $\displaystyle\sum_{k=2}^{a-4}c_{k}k=ta+(a-1)$, in the equation \ref{eq5.4}, we get\\
 $b=(\displaystyle \sum_{k=2}^{a-4}c_{k}-1)a+(ta+a-3)d$.\\
If $t\geq 2$ then $(\displaystyle \sum_{k=2}^{a-4}c_{k}-1)a+(ta+a-3)d> 2a+(a-3)d$. Since $b< 2a+(a-3)d$, which is a contradiction. Therefore $0\leq t\leq 2$.\\
If $t=0$ then $\displaystyle \sum_{k=2}^{a-4}c_{k}>2$ gives contradiction to the $b< 2a+(a-3)d$, hence $0<\displaystyle \sum_{k=2}^{a-4}c_{k}\leq 2$.
\end{enumerate}
\begin{enumerate}[(i)]
\item If any $c_{j}\neq 0$, for some $2\leq j \leq a-4$, then by the equation \ref{eq5.3} we get $b=(
k-2)d$, where $0\leq k-2\leq a-6$. $b\equiv (k-2)i (\mathrm{mod}\, a)$ which is a contradiction to $b\equiv (a-3)i (\mathrm{mod}\, a)$. 
\item If any $c_{j}\neq 0$ and $c_{\ell}\neq 0$ where $2\leq j,\ell\leq a-4$ such that $j\neq \ell$, then $b-a=(k_{1}+k_{2}-2)d$. We get $d|b-a$ which is a contradiction.
\end{enumerate}
Hence the element $a+b+2d\in\mathrm{Ap}(\Gamma_{e}(M),a)\setminus M$, can be written as $a+b+2d=(a+d)+(b+d)=(a+2d)+b$. Therefore the order of $a+b+2d$ is $2$.
\smallskip

Similarly we can prove if $b\equiv(a-2)i(\mathrm{mod}\, a)$ the element $2a+(a-3)d\in \mathrm{Ap}(\Gamma_{e}(M),a)\setminus M$, can be expressed as $2a+(a-3)d=a+k_{1}d+a+k_{2}d$ such that $k_{1}+k_{2}=a-3$. Therefore, the order of $2a+(a-3)d$ is $2$. \qed
\end{enumerate}
  
\begin{corollary} Let $\Gamma_{e}(M)(t_{k})$ be the number 
of elements of order $k$ ($1\leq k\leq 2$) in $\mathrm{Ap}(\Gamma_{e}(M),a)$. 
\begin{enumerate}
\item If $e=4$, then in both the cases $b\equiv 2i(\mathrm{mod}\,5)$ and $b\equiv 3i(\mathrm{mod}\, 5)$, 
$\Gamma_{4}(M)(t_{k})$ is given as follows:
\begin{table}[ht]
%\caption{No of elements in Ap\'{e}ry set of particular order} % title of Table
\centering % used for centering table
\begin{tabular}{|c|c|c|} % centered columns (4 columns)
\hline %inserts double horizontal lines
Value of $k$ & $1$ & $2$ \\ [0.5ex] % inserts table
%heading
\hline % inserts single horizontal line
$\Gamma_{4}(M)(t_{k})$  & $3$ & $1$ \\\hline % inserting body of the table
\end{tabular} % is used to refer this table in the text
\end{table}

\item If $e\geq 5$, then in both the cases $b\equiv (a-3)i (\mathrm{mod}\,5)$ and $b\equiv (a-2)i (\mathrm{mod}\, 5)$, 
$\Gamma_{e}(M)(t_{k})$ is given as follows:
\begin{table}[ht]
%\caption{No of elements in Ap\'{e}ry set of particular order} % title of Table
\centering % used for centering table
\begin{tabular}{|c|c|c|} % centered columns (4 columns)
\hline %inserts double horizontal lines
Value of $k$ & $1$ & $2$ \\ [0.5ex] % inserts table
%heading
\hline % inserts single horizontal line
$\Gamma_{e}(M)(t_{k})$  & $e-1$ & $1$ \\\hline % inserting body of the table
\end{tabular} % is used to refer this table in the text
\end{table}

\end{enumerate}
\end{corollary}
\proof Follows from Theorem \ref{orderam}.\qed 

\begin{corollary}\label{exptangent1} Let $I=\langle t^{a}\rangle$. The tangent cone $G_{\mathfrak{m}}(\Gamma_{e}(M))$ of $\Gamma_{e}(M)$ is a free $F(I)$-module. The 
following statements are true:
\begin{enumerate}
\item $ G_{\mathfrak{m}}(\Gamma_{4}(M))= \displaystyle\bigoplus(F(I)(-1))^{3}\bigoplus(F(I)(-2))^{1}$.
\item For $e\geq 5$, $ G_{\mathfrak{m}}(\Gamma_{e}(M))= \displaystyle\bigoplus(F(I)(-1))^{e-1}\bigoplus(F(I)(-2))^{1}.$
\end{enumerate} 
\end{corollary}

\proof Proof follows from the theorems \ref{tangentcone}, \ref{unique4}, \ref{aperytable4} and Lemma \ref{order4}.\qed
 
\begin{corollary} Tangent cones $G_{\mathfrak{m}}(\Gamma_{4}(M))$ and $G_{\mathfrak{m}}(\Gamma_{e}(M))$, $e\geq 5$ are Cohen-Macaulay.
\end{corollary}

\proof It easily follows from the fact that 
$G_{\mathfrak{m}}(\Gamma_{e}(M))$ is a free 
$F(I)$-module (see section 4 in \cite{cz2}). \qed

\begin{corollary}Let $HG_{\mathfrak{m}}(\Gamma_{e}(M))(x)$ be the Hilbert series of $G_{\mathfrak{m}}(\Gamma_{e}(M))$, 
then following statements are true: 
\begin{enumerate}
\item $HG_{\mathfrak{m}}(\Gamma_{4}(M))(x)=\displaystyle\left(3x+x^{2}\right)/(1-x)$.
\item For $e\geq 5$, 
$HG_{\mathfrak{m}}(\Gamma_{e}(M))(x)=\displaystyle\left((e-1)x+x^{2}\right)/(1-x).$
\end{enumerate}
\end{corollary}

\proof Follows from Corollay \ref{exptangent1}.\qed

\section{Tangent cone of the unbounded class $\mathfrak{S}_{(n,j)}, j=4,5$}
Let $e\geq 4$, $p\geq 0$ and $n=(e-1)+k(e-3)$, where $k\geq 2$. We have defined in \cite{mss3} the 
numerical semigroup $\mathfrak{S}_{(n,e,p)}$, generated by the 
integers $\{m_{0},\ldots,m_{e-1}\}$, where 
$m_{i}:=n^2+(e-2)n+p+i$, for $0\leq i\leq e-3 $ and $m_{e-2}:=n^2+(e-1)n+p+(e-3)$, 
$m_{e-1}:=n^2+(e-1)n+p+(e-2)$. This is formed by concatenation of two arithmetic 
sequences with common difference $1$. If $p=e-4$, then we write $\mathfrak{S}_{(n,e,e-4)}$ instead of $\mathfrak{S}_{(n,e)}$. Let  
$\mathcal{Q}_{(n,e,p)}\subset k[x_{0},\ldots,x_{e-1}]$ be the defining ideal of 
$\mathfrak{S}_{(n,e,p)}$.  
These classes are special in the sense that there is no  
upper bound on the minimal number of generators of the 
ideal $\mathcal{Q}_{(n,e,p)}$, for $e = 4, 5$.

\medskip

\begin{theorem}\label{Apery 4}
The \textit{Ap\'{e}ry set} of $\mathfrak{S}_{(n,4)}$ with respect to  $m_{0}=n^{2}+2n$ is
$\rm{Ap}(\mathfrak{S}_{(n,4)},m_{0})=\displaystyle\cup_{i=1}^{5} A_{i}\cup \{0\}$, where
\begin{align*}
A_{1}&=\{rm_{1}\mid 1\leq r\leq n\},\\
A_{2}&=\{rm_{2}\mid 1\leq r\leq n\},\\
A_{3}&=\{rm_{3}\mid 1\leq r\leq n-1\},\\
A_{4}&=\{rm_{1}+sm_{3}\mid 1\leq r\leq n-1,1\leq s\leq n-r\},\\
A_{5}&=\{rm_{2}+sm_{3}\mid 1\leq r\leq n-1,1\leq s\leq n-r\}.
\end{align*}
\end{theorem}
\proof See the theorem $3.2$\,\cite{mss5}. \qed
\medskip

\begin{theorem}\label{unique4}
Every element in $\rm{Ap}(\mathfrak{S}_{(n,4)},m_{0})$ has unique expression.
\end{theorem}
\proof Since $m_{1}$ is an element of minimal generating set for the semigroup $\mathfrak{S}_{(n,4)}$, so it has unique expression. Suppose $nm_{1}=c_{1}m_{1}+c_{2}m_{2}+c_{3}m_{3}$ (since $nm_{1}\in \rm{Ap}(\mathfrak{S}_{(n,4)},m_{0}) $ coefficient of $m_{0}$ is zero). If $0<c_{1}<n$, then $(n-c_{1})m_{1}=c_{2}m_{2}+c_{3}m_{3}$. But by induction $(n-c_{1})m_{1}$ has unique expression, we get a contradiction. If $c_{1}= 0$ then 
$$nm_{1} =c_{2}m_{2}+c_{3}m_{3} =c_{2}(m_{1}+n)+c_{3}(m_{1}+n+1).$$
Therefore $(n-(c_{1}+c_{2}))m_{1} = c_{2}n+c_{3}(n+1).$ Since R.H.S. is positive we have $c_{1}+c_{2}<n$, hence  $c_{2}n+c_{3}(n+1)<n^{2}+n$ and we get a contradiction. By similar approach we can show elements are in $A_{2}\cup A_{3}$ have unique expression.
\smallskip

Next we want to show that elements are in $A_{4}$ have unique expression. Fix $r$, $1\leq r\leq n$. It is enough to show that $rm_{1}+(n-r)m_{3}$ has unique expression. Suppose $rm_{1}+(n-r)m_{3}=c_{1}m_{1}+c_{2}m_{2}+c_{3}m_{3}$. Since $m_{2}=m_{1}+n$ and $m_{3}=m_{1}+n+1$, we get, 
$$(n-(c_{1}+c_{2}+c_{3}))m_{1}+(n-r-c_{3})(n+1)=c_{2}n.$$ 
Since R.H.S. is positive we have $c_{1}+c_{2}+c_{3}\leq n$. If $n=c_{1}+c_{2}+c_{3}$ then $(n-r-c_{3})(n+1)=c_{2}n $. Hence $(n+1)\mid c_{2}$ and we get $c_{2}=0$. Therefore $c_{3}=n-r$ and $c_{1}=r$. If $ c_{1}+c_{2}+c_{3}< n$ then, 
$$(n-(c_{1}+c_{2}+c_{3}))m_{1}=c_{2}n-(n-r-c_{3})(n+1)$$ 
and we have $c_{2}n-(n-r-c_{3})(n+1)\leq n^{2}-n<m_{1}$, therefore we get a contradiction. By a similar approach we can show elements are in $A_{5}$ have unique expression.  \qed
\medskip

\begin{lemma}\label{order4}
With the notion of the theorem \ref{Apery 4} for $k\geq 1$,
\begin{enumerate}
\item[(i)] $\mathrm{ord}(rm_{i}+km_{0})=r+k$, for all $1\leq r\leq n, i=1,2$
\item[(ii)] $\mathrm{ord}(rm_{3}+km_{0})=r+k$, for all $1\leq r\leq n-1$
\item[(iii)] $\mathrm{ord}((rm_{i}+sm_{3}+km_{0})=r+k+s$, for all $1\leq r\leq n-1,1\leq s\leq n-r, i=1,2$
\end{enumerate}
\end{lemma}
\proof   At first we show that $\mathrm{ord}(rm_{1}+sm_{3}+m_{0})=r+s+1$ for all $1\leq r\leq n,0\leq s\leq n-r$. Suppose for some $1\leq r\leq n,0\leq s\leq n-r$
\begin{equation}\label{eq1}
rm_{1}+sm_{3}+m_{0}=\sum_{i=0}^{3}l_{i}m_{i} 
\end{equation}
If $l_{0}>0$ then we have $rm_{1}+sm_{3}=(l_{0}-1)m_{0}+ \sum_{i=1}^{3}l_{i}m_{i}$ and by the theorem \ref{unique4}, $l_{0}=1,l_{1}=r,l_{2}=0,l_{3}=s$ and we are done.
\smallskip

If $l_{0}=0$ then rearranging the equation \ref{eq1}, we get 
\begin{equation}\label{eq2}
r+s(n+2)=(l_{1}+l_{2}+l_{3}-r-s-1)m_{0}+ l_{1}+l_{2}(n+1)+l_{3}(n+2)
\end{equation}

Again $1\leq r\leq n, 0\leq s\leq n-r$ imply $r+s(n+2)\leq r+(n-r)(n+2)\leq 1+(n-1)(n+2)<m_{0}$. Therefore from the equation \ref{eq2}, $l_{1}+l_{2}+l_{3}\leq r+s+1$.
\smallskip

Next we assume for $k>1$, 
\begin{equation}\label{eq3}
rm_{1}+sm_{3}+km_{0}=\sum_{i=0}^{3}c_{i}m_{i} 
\end{equation} 

If $c_{0}\geq k$ then $ rm_{1}+sm_{3}=(c_{0}-k)m_{0}+\sum_{i=1}^{3}c_{i}m_{i}$ and by the theorem \ref{unique4}, $c_{0}=k,l_{1}=r,l_{2}=0,l_{3}=s$ and we are done.
\smallskip

If $1\leq c_{0}\leq k$ then $ rm_{1}+sm_{3}+(k-c_{0})m_{0}=\sum_{i=1}^{3}c_{i}m_{i}$ and by induction we have $\sum_{i=1}^{3}c_{i}\leq r+s+(k-c_{0})$. If $c_{0}=0$, then we have 
\begin{equation}\label{eq3a}
r+s(n+2)=(c_{1}+c_{2}+c_{3}-r-s-k)m_{0}+ c_{1}+c_{2}(n+1)+c_{3}(n+2)
\end{equation}
Again $1\leq r\leq n, 0\leq s\leq n-r$ implies $r+s(n+2)<m_{0}$. Therefore from the equation \ref{eq3a}, $l_{1}+l_{2}+l_{3}\leq r+s+k$. which completes the proof that $\mathrm{ord}(rm_{1}+sm_{3}+km_{0})=r+s+k$, for $k\leq 1$ and $1\leq r\leq n,0\leq s\leq n-r$.
\smallskip

Similarly we can prove that $\mathrm{ord}(rm_{2}+sm_{3}+km_{0})=r+s+k$, for $k\geq 1$ and $1\leq r\leq n,0\leq s\leq n-r$.

\medskip

Next we prove that $\mathrm{ord}(rm_{3}+km_{0})=r+k$, for all $1\leq r\leq n-1$ and $k\geq 1$. Again at first we show that $\mathrm{ord}(rm_{3}+m_{0})=r+k$, for all $1\leq r\leq n-1$. Suppose for some $1\leq r\leq n-1$
\begin{equation}\label{eq4}
rm_{3}+m_{0}=\sum_{i=0}^{3}l_{i}m_{i} 
\end{equation}
If $l_{0}>0$ then we have $rm_{3}=(l_{0}-1)m_{0}+ \sum_{i=1}^{3}l_{i}m_{i}$ and by the theorem \ref{unique4}, $l_{0}=1,l_{1}=0,l_{2}=0,l_{3}=r$ and we are done.
\smallskip

If $l_{0}=0$ then rearranging the equation \ref{eq4}, we get 
\begin{equation}\label{eq5}
r(n+2)=(l_{1}+l_{2}+l_{3}-r-1)m_{0}+ l_{1}+l_{2}(n+1)+l_{3}(n+2)
\end{equation}

 Again $1\leq r\leq n-1$ implies $r(n+2)\leq (n-1)(n+2)<m_{0}$. Therefore from the equation \ref{eq5}, $l_{1}+l_{2}+l_{3}=r+s+1$.
\smallskip

Suppose $k>1$, 
\begin{equation}\label{eq6}
rm_{3}+km_{0}=\sum_{i=0}^{3}c_{i}m_{i} 
\end{equation} 

If $c_{0}\geq k$ then $ rm_{3}=(c_{0}-k)m_{0}+\sum_{i=1}^{3}c_{i}m_{i}$ and by the theorem \ref{unique4}, $c_{0}=k,l_{1}=r,l_{2}=0,l_{3}=s$ and we are done.
\smallskip

If $1\leq c_{0}\leq k$ then $ rm_{3}+(k-c_{0})m_{0}=\sum_{i=1}^{3}c_{i}m_{i}$ and by induction we have $\sum_{i=1}^{3}c_{i}\leq r+(k-c_{0})$ and we are done. 
\smallskip

If $c_{0}=0$, then we have
\begin{equation}\label{eq7}
r(n+2)=(c_{1}+c_{2}+c_{3}-r-k)m_{0}+ c_{1}+c_{2}(n+1)+c_{3}(n+2)
\end{equation}
 Since $r(n+2)<m_{0}$ for $1\leq r\leq n-1$, from the equation \ref{eq7} we get $c_{1}+c_{2}+c_{3}<r+k $, which completes the proof that $\mathrm{ord}(rm_{3}+km_{0})=r+s+k$, for $k\leq 1$ and $1\leq r\leq n-1$.
\qed
\medskip

\begin{theorem}\label{aperytable4}
The Ap\'{e}ry table $\mathrm{AT}(\mathfrak{S}_{(n,4)},m_{0})_{(n+1)\times m_{0}}$ of  $\mathfrak{S}_{(n,4)}$ 
is of the form 
$\mathrm{AT}(\mathfrak{S}_{(n,4)},m_{0})=\begin{pmatrix}
\mathcal{A}_{0}& \mathcal{A}_{1}& \mathcal{A}_{2}& \mathcal{A}_{3}& \mathcal{A}_{4}&\mathcal{A}_{5}
\end{pmatrix}$, where 
\begin{itemize}
\item $\mathcal{A}_{0}=(a_{i1}^{(0)})_{(n+1)\times 1},\,\, a_{i1}^{(0)}=(i-1)m_{0},\, 1\leq i\leq n+1$
\item $\mathcal{A}_{1}=(a_{ij}^{(1)})_{(n+1)\times n}$, where\,\, $a_{ij}^{(1)}=jm_{1}+\gamma_{ij}^{(1)}m_{0}$ and
\begin{align*}
\gamma_{ij}^{(1)}&=(i-j-1) \,\, \rm{if}\, i>j+1\\
&=0 \,\, \rm{if}\, i\leq j+1
\end{align*}
\item $\mathcal{A}_{2}=(a_{ij}^{(2)})_{(n+1)\times n}$, where\,\, $a_{ij}^{(2)}=jm_{2}+\gamma_{ij}^{(2)}m_{0}$ and
\begin{align*}
\gamma_{ij}^{(2)}&=(i-j-1) \,\, \rm{if}\, i>j+1\\
&=0 \,\, \rm{if}\, i\leq j+1
\end{align*}
\item $\mathcal{A}_{3}=(a_{ij}^{(3)})_{(n+1)\times (n-1)}$, where\,\, $a_{ij}^{(1)}=jm_{3}+\gamma_{ij}^{(3)}m_{0}$ and
\begin{align*}
\gamma_{ij}^{(3)}&=(i-j-1) \,\, \rm{if}\, i>j+1\\
&=0 \,\, \rm{if}\, i\leq j+1
\end{align*}
\item $\mathcal{A}_{4}=\begin{pmatrix}
\Omega_{14}&\cdots & \Omega_{(n-1)4}
\end{pmatrix}$ Where 
$\Omega_{r4}=(\omega_{ij}^{(r4)})_{(n+1)\times (n-r)}, 1\leq r\leq n-1$ and  for  $1\leq i\leq n+1, 1\leq j\leq n-r$, we have $\omega_{ij}^{(r4)}=rm_{1}+jm_{3}+\alpha_{ij}^{(r4)}m_{0}$,
\begin{align*}
\alpha_{ij}^{(r4)}&=i-j-1-r,\quad \rm{if}\, i-j-1-r>0\\
&=0 \quad \rm{if}\, i-j-1-r\leq 0
\end{align*}
\item $\mathcal{A}_{5}=\begin{pmatrix}
\Omega_{15}&\cdots & \Omega_{(n-1)5}
\end{pmatrix}$ Where 
$\Omega_{r5}=(\omega_{ij}^{(r5)})_{(n+1)\times (n-r)}, 1\leq r\leq n-1$ and  for  $1\leq i\leq n+1, 1\leq j\leq n-r$, we have $\omega_{ij}^{(r5)}=rm_{2}+jm_{3}+\alpha_{ij}^{(r5)}m_{0}$,
\begin{align*}
\alpha_{ij}^{(r5)}&=i-j-1-r,\quad \rm{if}\, i-j-1-r>0\\
&=0 \quad \rm{if}\, i-j-1-r\leq 0
\end{align*}
\end{itemize}
\end{theorem}
\proof Follows from Theorem \ref{unique4} and Lemma \ref{order4}.\qed
\medskip

\begin{corollary}\label{table4}
Let $\mathfrak{S}_{(n,4)}(t_{k})$ be the number of 
elements of particular order $k$, $1\leq k\leq n$, 
in the Ap\'{e}ry set $\rm{Ap}(\mathfrak{S}_{(n,4)},m_{0})$. Then $\mathfrak{S}_{(n,4)}(t_{k})=2k+1$, for $1\leq k\leq n$.
\end{corollary}

\proof We give a table which follows from the Theorems  \ref{unique4} and \ref{aperytable4}:
\begin{table}[ht]
 \centering 
\begin{tabular}{|c|c|c|c|c|c|c|c|}
\hline
Values of $k$ & $1$ & $2$ & $3$& $\cdots$ & $r$ & $\cdots$ & $n$\\ [0.5ex] 
\hline 
$A_{1}$ & $1$ & $1$ & $1$ &$\cdots$ & $1$ &$\cdots$ & $1$\\\hline 
$A_{2}$ & $1$ & $1$ & $1$ &$\cdots$& $1$ & $\cdots$& $1$\\\hline
$A_{3}$ & $1$ & $1$ & $1$ &$\cdots$& $1$ &$\cdots$& $1$\\\hline
$A_{4}$ & $0$ & $1$ & $2$ &$\cdots$ & $r-1$ &$\cdots$ &$n-1$ \\\hline
$A_{5}$ & $0$ & $1$ & $2$ &$\cdots$ &$r-1$ &$\cdots$ &$n-1$\\ \hline
$\mathfrak{S}_{(n,4)}(t_{k})$ & $3$ & $5$ & $7$ &$\cdots$ & $2r+1$ & $\cdots$ & $2n+1$\\ \hline 
\end{tabular} \qed

\end{table} 
\medskip

\begin{corollary}\label{exptangent2} Let $I=\langle t^{n^{2}+2n}\rangle$. The tangent cone $G_{\mathfrak{m}}(\mathfrak{S}_{(n,4)})$ of $\mathfrak{S}_{(n,4)}$ is a free 
$F(I)$-module. Moreover,  
$ G_{\mathfrak{m}}(\mathfrak{S}_{(n,4)})= \displaystyle\bigoplus_{k=1}^{n}(F(I)(-k))^{2k+1}$.  
\end{corollary}  

\proof Proof follows from Theorems \ref{tangentcone}, \ref{unique4}, \ref{aperytable4} and Lemma \ref{order4}.\qed
\medskip

\begin{corollary} The tangent cone $G_{\mathfrak{m}}(\mathfrak{S}_{(n,4)})$ is Cohen-Macaulay.
\end{corollary}

\proof It easily follows                                    from the fact that  $G_{\mathfrak{m}}(\mathfrak{S}_{(n,4)})$ is a free $F(I)$-module (see section 4 in \cite{cz2}). \qed
\medskip

\begin{corollary} Let $HG_{\mathfrak{m}}(\mathfrak{S}_{(n,4)})(x)$ be the Hilbert series of $G_{\mathfrak{m}}(\mathfrak{S}_{(n,4)})$. Then $HG_{\mathfrak{m}}(\mathfrak{S}_{(n,4)})(x)=\displaystyle\left(\sum_{k=1}^{n} (2k+1)x^{k}\right)/(1-x)$.
\end{corollary}
\proof Follows from Corollay \ref{exptangent2}.\qed
\medskip

\begin{example}
 For $n=5$, $\mathfrak{S}_{(5,4)}=\langle 35,36,41,42 \rangle$. The Ap\'{e}ry set $\rm{Ap}(\mathfrak{S}_{(5,4)},m_{0})=\displaystyle\cup_{i=1}^{5} A_{i}\cup \{0\}$, where
\begin{align*}
A_{1}&=\{36,72,108,144,180\},\\
A_{2}&=\{41,82,123,164,205\},\\
A_{3}&=\{42,84,126,168\},\\
A_{4}&=\{78,114,120,150,156,162,186,192,198,204\},\\
A_{5}&=\{83,124,125,165,166,167,206,207,208,209\}.
\end{align*}
The Ap\'{e}ry table $\mathrm{AT}(\mathfrak{S}_{(5,4)},35)_{6\times 35}$ of  $\mathfrak{S}_{(5,4)}$ is of the form 
$$\mathrm{AT}(\mathfrak{S}_{(5,4)},35)=\begin{pmatrix}
\mathcal{A}_{0}& \mathcal{A}_{1}& \mathcal{A}_{2}& \mathcal{A}_{3}& \mathcal{A}_{4}&\mathcal{A}_{5}
\end{pmatrix},$$ 
where 

\begin{tabular}{cc}   % top level tables, with 2 columns
%$(i)$ & $(ii)$\\  

(i)\,  $\mathcal{A}_{0}= \begin{tabular}{|c|}
\hline
 0  \\ 
 \hline
 35\\
 \hline 
 70\\
 \hline
 105\\
 \hline
 140\\ 
 \hline
 175\\
 \hline   
\end{tabular}
$
\quad \quad \quad \quad \quad \quad \quad \quad  \quad \quad \quad (ii)\, $\mathcal{A}_{1}=\begin{tabular}{|c|c|c|c|c|}
\hline
36& 72 & 108 & 144 & 180\\
\hline
36& 72 & 108 & 144 & 180\\
\hline
71& 72 & 108 & 144 & 180\\
\hline
106& 107 & 108 & 144 & 180\\
\hline
141& 142 & 143 & 144 & 180\\
\hline
176& 177 & 178 & 179 & 180\\
\hline   
\end{tabular}
$
\end{tabular}

\bigskip

\begin{tabular}{cc}
(iii) $\mathcal{A}_{2}=\begin{tabular}{|c|c|c|c|c|}
\hline
41& 82 & 123 & 144 & 205\\
\hline
41& 82 & 123 & 144 & 205\\
\hline
76& 82 & 123 & 144 & 205\\
\hline
111& 117 & 123 & 144 & 205\\
\hline
146& 152 & 158 & 144 & 205\\
\hline
181& 187 & 193 & 179 & 205\\
\hline   
\end{tabular}
$
\quad (iv)\,  $\mathcal{A}_{3}=\begin{tabular}{|c|c|c|c|c|}
\hline
42& 84 & 126 & 168 \\
\hline
42& 84 & 126 & 168\\
\hline
77& 84 & 126 & 168\\
\hline
112& 119 & 126 & 168\\
\hline
147& 154 & 161 & 168\\
\hline
182& 189 & 196 & 203\\
\hline   
\end{tabular}
$
\end{tabular}
\bigskip

(v)\,  $\mathcal{A}_{4}= \begin{pmatrix}
\Omega_{14}&\cdots & \Omega_{44}
\end{pmatrix}$, 
\smallskip

such that  

\begin{tabular}{cc}
(a)\,  $\Omega_{14}=\begin{tabular}{|c|c|c|c|c|}
\hline
78& 120 & 162 & 206 \\
\hline
78& 120 & 162 & 206\\
\hline
78& 120 & 162 & 206\\
\hline
113& 120 & 162 & 206\\
\hline
148& 155 & 162 & 206\\
\hline
183& 190 & 197 & 206\\
\hline   
\end{tabular}
$
\quad \quad \quad (b)\,  $\Omega_{24}=\begin{tabular}{|c|c|c|c|c|}
\hline
114& 156 & 198  \\
\hline
114& 156 & 198 \\
\hline
114& 156 & 198 \\
\hline
114& 156 & 198 \\
\hline
149& 156 & 198 \\
\hline
184& 191 & 198 \\
\hline   
\end{tabular}
$
\end{tabular}
\bigskip

\begin{tabular}{cc}
 (c)\, $\Omega_{34}=\begin{tabular}{|c|c|}
\hline
150& 192 \\
\hline
150& 192 \\
\hline
150& 192  \\
\hline
150& 192 \\
\hline
150& 192  \\
\hline
185& 192  \\
\hline   
\end{tabular}
$
\quad \quad \quad \quad \quad \quad \quad \quad (d)\,  $\Omega_{44}=\begin{tabular}{|c|}
\hline
186\\
\hline
186\\
\hline
186\\
\hline
186\\
\hline
186\\
\hline
186\\
\hline   
\end{tabular}
$
 \end{tabular}
\end{example}

\begin{theorem}\label{Apery 5}
The \textit{Ap\'{e}ry set} of $\mathfrak{S}_{(n,5)}$ with respect to  $m_{0}=n^{2}+3n+1$ is
$\mathrm{Ap}(\mathfrak{S}_{(n,5)},m_{0})=\displaystyle\cup_{i=1}^{11} A_{i}$, where
\begin{align*}
A_{1}& =\{0,m_{1}\},\\
A_{2}&=\{rm_{2}\mid 1\leq r\leq \dfrac{n}{2}\},\\
A_{3}&=\{rm_{3}\mid 1\leq r\leq n\},\\
A_{4}&=\{rm_{4}\mid 1\leq r\leq n\}\\
A_{5}&=\{m_{1}+r m_{2}\mid 1\leq r\leq \dfrac{n}{2}\},\\
A_{6}&=\{m_{3}+r m_{2}\mid 1\leq r\leq \dfrac{n}{2}\}\\
A_{7}&=\{rm_{2}+2s m_{4}\mid 1\leq s\leq \dfrac{n}{2}-1,1\leq r \leq \dfrac{n}{2}-s \},\\
A_{8}&=\{rm_{2}+(2s-1)m_{4}\mid 1\leq s\leq \dfrac{n}{2},1\leq r \leq \dfrac{n}{2}+1-s \},\\
A_{9}&=\{km_{3}+(n-k-r+1)m_{4}\mid 1\leq k\leq n-1,1\leq r \leq n-k \},\\
A_{10}&=\{rm_{2}+m_{3}+2s m_{4}\mid 1\leq s\leq \dfrac{n}{2}-1,1\leq r \leq \dfrac{n}{2}-s \},\\
A_{11}&=\{rm_{2}+m_{3}+(2s-1)m_{4}\mid 1\leq s\leq \dfrac{n}{2}-1,1\leq r \leq \dfrac{n}{2}+1-s \}.
\end{align*}
\end{theorem}
\proof See Theorem 4.2 \cite{mss5}. \qed
\medskip

\begin{theorem}\label{unique5}
The set $\rm{Ap}(\mathfrak{S}_{(n,5)},m_{0})$ is homogeneous.
\end{theorem}
\proof By a similar way as in Theorem \ref{unique4}, we can show that each element of $\cup_{i=1}^{9}A_{i}$ has unique expression. Next we want to show that order of an element $rm_{2}+m_{3}+2s m_{4}$, $ 1\leq s\leq \dfrac{n}{2}-1$, 
$1\leq r \leq \dfrac{n}{2}-s$ in $A_{10}$ is $r+2s+1$. 

Suppose  $rm_{2}+m_{3}+2sm_{4}=c_{0}m_{0}+c_{1}m_{1}+c_{2}m_{2}+c_{3}m_{3}+c_{4}m_{4}$, then 
\begin{align*}
&(r+2s+1-c_{0}-c_{1}-c_{2}-c_{3}-c_{4})m_{0}\\
&=c_{1}+2c_{2}+(c_{3}-1)(n+2)+(c_{4}-2s)(n+3)-2r.
\end{align*}

\noindent\textbf{Case 1.} If $c_{0}+c_{1}+c_{2}+c_{3}+c_{4}< r+2s+1$. Since $r+2s+1\leq n-1$, we have
\begin{align*}
R.H.S &= c_{1}+2c_{2}+(c_{3}-1)(n+2)+(c_{4}-2s)(n+3)-2r\\
&= c_{1}+2c_{2}+c_{3}(n+2)+c_{4}(n+3)-n-2sn-2(r+3s+1)\\
&\leq 2c_{0}+2c_{1}+2c_{2}+2c_{3}+c_{3}n+c_{4}(n+1)-n-2sn-2(r+3s+1)\\
&< 2(r+2s+1)+c_{3}n+c_{4}(n+1)-n-2sn-2(r+3s+1)\\
&=(c_{3}+c_{4})n+c_{4}-n-2sn-2s\\
&<(r+2s+1)(n+1)\leq (n-1)(n+1)=n^{2}-1< n^{2}+3n+1=m_{0}.
\end{align*}
This gives a contradiction. 
\medskip

\noindent\textbf{Case 2.} If $c_{0}+c_{1}+c_{2}+c_{3}+c_{4}> r+2s+1$. We have $r+3s+1\leq 3(\dfrac{n}{2}-1)$, 
\begin{align*}
&R.H.S=c_{1}+2c_{2}+(c_{3}-1)(n+2)+(c_{4}-2s)(n+3)-2r\\
&=(c_{1}+c_{2}+c_{3}+c_{4})+(c_{2}+c_{3}+2c_{4})+n(c_{3}-1)+n(c_{4}-2s)-6s-2r-2\\
&\geq n(c_{3}-1)+n(c_{4}-2s)-6s-2r-2\\
&\geq -n-2sn-2r-6s-2\geq -n-(\dfrac{n}{2}-2)n-6(\dfrac{n}{2}-1)\\
&=-\dfrac{n^{2}}{2}-2n+6>-n^{2}-3n-1=-m_{0}.
\end{align*}

\noindent We get a contradiction. Therefore $c_{0}+c_{1}+c_{2}+c_{3}+c_{4}= r+2s+1$ and the order of an element $rm_{2}+m_{3}+2s m_{4}$, $ 1\leq s\leq \dfrac{n}{2}-1,1\leq r \leq \dfrac{n}{2}-s$ in $A_{10}$ is $r+2s+1$.  
\medskip

By a similar approach, we can show that the order of an element $rm_{2}+m_{3}+(2s-1) m_{4}$, $ 1\leq s\leq \dfrac{n}{2}-1,1\leq r \leq \dfrac{n}{2}+1-s$ in $A_{11}$ is $r+2s$. 
\medskip

\begin{lemma}\label{order5}
With the notion of Theorem \ref{Apery 5}, for $k\geq 1$,
\begin{enumerate}
\item[(i)] ord$(m_{1}+km_{0})=k+1$;
\item[(ii)] ord$(rm_{2}+km_{0})=r+k$, for all $1\leq r\leq \dfrac{n}{2}$;
\item[(iii)] ord$(rm_{i}+km_{0})=r+k$, for all $1\leq r\leq n, i=3,4$;
\item[(iv)] ord$(m_{i}+rm_{3}+km_{0})=r+k+1$, for all $1\leq r\leq \dfrac{n}{2}, i=1,3$;
\item[(v)] ord$((rm_{2}+sm_{3}+km_{0})=r+k+s$, for all $1\leq r\leq n-1,1\leq s\leq n-r$.
\end{enumerate}
\end{lemma}
\proof The proof is similar as in \ref{order4}.\qed
\medskip

\begin{theorem}\label{AT5}
The Ap\'{e}ry table $\mathrm{AT}(\mathfrak{S}_{(n,5)},m_{0})_{(n+1)\times m_{0}}$ of  $\mathfrak{S}_{(n,5)}$ is of the form 
$\mathrm{AT}(\mathfrak{S}_{(n,4)},m_{0})=\begin{pmatrix}
\mathcal{B}_{0}& \cdots & \mathcal{B}_{11}
\end{pmatrix}$, where 
\begin{itemize}
\item $\mathcal{B}_{1}=(b_{i1}^{(1)})_{(n+1)\times 2}$, 
such that\\

$b_{i1}^{(1)} = 
\begin{cases}
(i-1)m_{0}, & \quad \rm{if} \quad  1\leq i\leq n+1, \, j=1;\\
m_{1}, & \quad \rm{if} \quad  i=1,2,\quad j=2;\\
m_{1}+(i-2)m_{0}, & \quad \rm{if} \quad  3\leq i\leq n+1,\quad j=2.
\end{cases}$

\item $\mathcal{B}_{2}=(b_{ij}^{(2)})_{(n+1)\times \frac{n}{2}}$, such that \, $b_{ij}^{(2)}=jm_{2}+\gamma_{ij}^{(2)}m_{0}$ \, and \\
$\gamma_{ij}^{(2)} = 
\begin{cases}
(i-j-1), & \quad \rm{if}\quad i>j+1;\\
0, & \quad \rm{if}\quad i\leq j+1.
\end{cases}$

\item For $k=3,4$, we have $\mathcal{B}_{k}=(b_{ij}^{(k)})_{(n+1)\times n}$\, , such that\,\, $b_{ij}^{(k)}=jm_{k}+\gamma_{ij}m_{0}$ \,\, and\\
$\gamma_{ij} = 
\begin{cases}
(i-j-1), & \quad \rm{if} \quad i>j+1;\\
0, & \quad \rm{if}\quad i\leq j+1.
\end{cases}$

\item For $k=5,6$, we have $\mathcal{B}_{k}=(b_{ij}^{(k)})_{(n+1)\times \frac{n}{2}}$, such that\,\, $b_{ij}^{(k)}=m_{k-3}+jm_{2}+\alpha_{ij}m_{0}$ \,\, and\\ 
$\alpha_{ij}= \begin{cases}
(i-j-2), & \quad \rm{if} \quad i>j+2;\\
0, & \quad \rm{if} \quad i\leq j+2.
\end{cases}
$

\item $\mathcal{B}_{7}=\begin{pmatrix}
Q_{17}&\cdots & Q_{(\frac{n}{2}-1)7}
\end{pmatrix}$, such that 
$Q_{s7}=(q_{ij}^{(s7)})_{(n+1)\times (\frac{n}{2}-s)}$, 
$q_{ij}^{(s7)}=2sm_{4}+jm_{2}+\delta_{ij}^{(s7)}m_{0} $ and\\
$\delta_{ij}^{(s7)}= 
\begin{cases}
(i-j-2s-1), & \quad \rm{if} \quad i>j+2s+1;\\
0, & \quad \rm{if} \quad i\leq j+2s+1.
\end{cases}$

\item $\mathcal{B}_{8}=\begin{pmatrix}
Q_{18}&\cdots & Q_{(\frac{n}{2})8}
\end{pmatrix}$, such that 
$Q_{s8}=(q_{ij}^{(s8)})_{(n+1)\times (\frac{n}{2}+1-s)}$, $q_{ij}^{(s8)}=(2s-1)m_{4}+jm_{2}+\delta_{ij}^{(s8)}m_{0}$ 
\, and \\ 
$\delta_{ij}^{(s8)}= 
\begin{cases}
(i-j-2s), & \quad \rm{if} \quad i>j+2s;\\
0, & \quad \rm{if} \quad i\leq j+2s.
\end{cases}$

\item $\mathcal{B}_{9}=\begin{pmatrix}
Q_{19}&\cdots & Q_{(n-1)9}
\end{pmatrix}$, such that $Q_{s9}=(q_{ij}^{(s9)})_{(n+1)\times (n-s)}$, $q_{ij}^{(s7)}=(n-s-j+1)m_{4}+jm_{3}+\delta_{ij}^{(s9)}m_{0}$ \, and\\
$\delta_{ij}^{(s9)} = 
\begin{cases}
(i+s-n-2), & \quad \rm{if} \quad i+s>n+2;\\
0, & \quad \rm{if} \quad i+s\leq n+2.
\end{cases}$

\item $\mathcal{B}_{10}=\begin{pmatrix}
Q_{110}&\cdots & Q_{(\frac{n}{2}-1)10}
\end{pmatrix}$, such that $Q_{s10}=(q_{ij}^{(s10)})_{(n+1)\times (\frac{n}{2}-s)}$, \, $q_{ij}^{(s10)}=2sm_{4}+jm_{2}+m_{3}+\delta_{ij}^{(s10)}m_{0}$ \, and\\
$\delta_{ij}^{(s10)} = 
\begin{cases}
(i-j-2s-2), & \quad \rm{if}\quad i>j+2s+2;\\
0, & \quad \rm{if}\quad i\leq j+2s+2.
\end{cases}$

\item $\mathcal{B}_{11}=\begin{pmatrix}
Q_{111}&\cdots & Q_{(\frac{n}{2}-1)11}
\end{pmatrix}$, such that 
$Q_{s11}=(q_{ij}^{(s11)})_{(n+1)\times (\frac{n}{2}+1-s)}$, 
\, 
$q_{ij}^{(s11)}=(2s-1)m_{4}+jm_{2}+m_{3}+\delta_{ij}^{(s11)}m_{0} $ and\\
$\delta_{ij}^{(s11)} = 
\begin{cases}
(i-j-2s-1), & \quad \rm{if}\quad i>j+2s+1;\\
0 & \quad \rm{if}\quad i\leq j+2s+1.
\end{cases}$

\end{itemize}
\end{theorem}

\proof Follows from Theorem \ref{unique5}. \qed

\begin{corollary}
Let $\mathfrak{S}_{(n,5)}(t_{k})$ be the number of elements of particular order $k$, $1\leq k\leq n$,  in the Ap\'{e}ry set $\rm{Ap}(\mathfrak{S}_{(n,5)},m_{0})$. Then
\begin{itemize}
\item $\mathfrak{S}_{(n,5)}(t_{1})=4$.

\item  $\mathfrak{S}_{(n,5)}(t_{2})=11$.

\item $\mathfrak{S}_{(n,5)}(t_{2m}) = 
\begin{cases}
6m+2, & \quad \rm{if} \quad 4\leq 2m\leq \dfrac{n}{2};\\
2n-2m+8, & \quad \rm{if} \quad \dfrac{n}{2}+1\leq 2m\leq n.
\end{cases}$

\item $\mathfrak{S}_{(n,5)}(t_{2m+1})= 
\begin{cases}
6m+5, & \quad \rm{if} \quad 3\leq 2m+1\leq \dfrac{n}{2};\\
2n-2m+7, & \quad \rm{if} \quad \dfrac{n}{2}+1\leq 2m+1\leq n-1.
\end{cases}$

\end{itemize}
\end{corollary}
\proof Follows from Theorem \ref{AT5}.\qed
\medskip

\begin{corollary}\label{exptangent3} Let $I=\langle t^{n^{2}+3n}\rangle$. The tangent cone $G_{\mathfrak{m}}(\mathfrak{S}_{(n,5)})$ of $\mathfrak{S}_{(n,5)}$ is a free 
$F(I)$-module. Moreover 
\begin{align*}
G_{\mathfrak{m}}(\mathfrak{S}_{(n,5)})&= \displaystyle \bigoplus(F(I)(-1))^4\bigoplus (F(I)(-2))^{11}\bigoplus_{m=1}^{n/4-1/2}(F(I)(-2m-1))^{6m+5}\\
&\bigoplus_{m=n/4}^{n/2-1}(F(I)(-2m-1))^{2n-2m+7}\bigoplus_{m=2}^{n/4}(F(I)(-2m))^{6m+2}\\
&\bigoplus_{m=n/4+1/2}^{n/2}(F(I)(-2m))^{2n-2m+8}.
\end{align*}
\end{corollary}

\proof Follows from the theorems \ref{tangentcone}, \ref{unique4}, \ref{aperytable4} and Lemma \ref{order4}.\qed
\medskip

\begin{corollary} Tangent cone $G_{\mathfrak{m}}(\mathfrak{S}_{(n,5)})$ is Cohen-Macaulay.
\end{corollary}

\proof It follows from the fact that  $G_{\mathfrak{m}}(\mathfrak{S}_{(n,5)})$ is a free $F(I)$-module (see section 4 in \cite{cz2}). \qed
\medskip

\begin{corollary} Let $HG_{\mathfrak{m}}(\mathfrak{S}_{(n,5)})(x)$ be the Hilbert series of $G_{\mathfrak{m}}(\mathfrak{S}_{(n,5)})$. Then 
\begin{align*}
HG_{\mathfrak{m}}(\mathfrak{S}_{(n,5)})(x)&=\displaystyle (4x+11x^{2})+\sum_{k=1}^{n/4-1/2} (6k+5)x^{2k+1}+\sum_{k=n/4}^{n/2-1}(2n-2k+7)x^{2k+1}\\
&+\sum_{k=2}^{n/4} (6k+2)x^{2k}+\sum_{k=n/4+1/2}^{n/2} (2n-2k+8)x^{2k})/(1-x).
\end{align*}
\end{corollary}
\proof Follows from Corollay \ref{exptangent3}.\qed

\bibliographystyle{amsalpha}

\end{document}